\documentclass[reqno, 12pt]{amsart}
\usepackage{amsfonts,amssymb}
\usepackage{enumerate}
\usepackage{verbatim}
\usepackage{mathrsfs}

\newcommand{\R}{\mathbb{R}}

\newcommand{\Z}{\mathbb{Z}}

\newtheorem{thm}{Theorem}
\newtheorem{lem}{Lemma}
\newtheorem{prop}{Proposition}

\theoremstyle{definition}

\theoremstyle{remark}

\numberwithin{equation}{section}

\author{Jing-Jing Huang}
\address{
Jing-Jing Huang: Department of Mathematics and Statistics, University of Nevada, Reno,
1664 N. Virginia St., Reno, NV 89557}
\email{jingjingh@unr.edu}
\dedicatory{}
\thanks{Research is supported by the UNR VPRI startup grant 1201-121-2479}

\begin{document}

\title
[Integral points close to a space curve]
{Integral points close to a  space curve}

\begin{abstract}
We establish sharp lower and upper bounds for the number of integral points near dilations of a space curve with nowhere vanishing torsion. 
\end{abstract}
\maketitle

\section{Integral points on a curve} \label{s1}

A classical theorem of Jarn\'{i}k \cite{J} states that a strictly convex arc in $\R^2$ of length ${q}\ge1$ contains at most $\ll q^{\frac23}$ integral points, and that the exponent $\frac23$ is best possible.

However, the above bound is susceptible to improvement if we start with a fixed strictly convex arc $\Gamma:y=f(x)$ and consider the number of integral points on the dilation $q\Gamma$ by a factor $q\ge1$. Indeed, Swinnerton-Dyer \cite{SD} showed that  for a $C^3$ strictly convex arc $\Gamma$ and for any $\varepsilon>0$
$$
\#(q\Gamma\cap\mathbb{Z}^2)\ll_{\Gamma,\varepsilon} q^{\frac35+\varepsilon}.
$$
Later on, W.M. Schmidt \cite{S} obtained a uniform version (with respect to $\Gamma$) of Swinnerton-Dyer's theorem and generalized it to hypersurfaces. In a landmark paper where an ingenious determinant method was developed, Bombieri and Pila \cite{BP} proved the optimal bound 
$$
\#(q\Gamma\cap\mathbb{Z}^2)\ll_{\Gamma,\varepsilon} q^{\frac12+\varepsilon}
$$
for strictly convex $f\in C^\infty[0,1]$, which confirms a conjecture of P. Sarnak. In a subsequent paper, Pila \cite{Pi} obtained the same bound under the weaker assumption  that $f\in C^{104}[0,1]$ and the determinant 
\[f''
\left|
\begin{array}{ccc}
f'''& 3f''& 0\\
f^{(4)}& 4f'''&6f''\\
f^{(5)}& 5f^{(4)}& 20f'''
\end{array}
\right|
\]
 is nowhere vanishing.
 
 Better results are available for algebraic curves of degree at least 3. Assuming $\Gamma$ is a subset of an irreducible algebraic curve of degree $d$ inside a square of side $q$, Bombieri and Pila showed, in the same paper cited above, that the number of lattice points on $\Gamma$ is 
 $$
 \ll_{d,\varepsilon} q^{1/d+\varepsilon}.
 $$
Similar bounds are also known for the problem of counting rational points on $\Gamma$, or equivalently counting integral points on the corresponding  projective curve \cite{EV,HB,W}.
 
 \section{Integral points close to a curve}
 
 Sometimes, we are not only interested in integral/rational points on a curve, but also those lying close to the curve, or equivalently those in a very thin neighborhood of the curve. 
The question of estimating the number of lattice points near dilations of the  unit circle $x^2+y^2=1$ is of course closely tied with the celebrated circle problem of Gauss, which asks for the best possible error term when approximating the number of lattice points inside  a circle centered at the origin with radius $r$ by its area $\pi r^2$. In general, one may investigate the number of lattice points near dilations of a reasonably smooth planar curve. For results on this rather difficult problem and applications to questions such as gaps between squarefull numbers, see \cite{hux2, T} and the references therein.  Swinnerton-Dyer's method \cite{SD} still plays an important role in this general setting, but the proofs are much more technical and even the results are too formidable to be reproduced here. However, for the standard parabola $\{(x,x^2), x\in[0,1]\}$, we have established a better bound by incorporating some multiplicative number theory pertaining to the arithmetic nature of the problem \cite{hua4}.

One may as well study this problem for proper submanifolds of the Euclidean space $\R^n$. In general, within a fixed ambient space, the problem becomes more difficult as the codimension of the underlying manifold increases. In fact, Lettington \cite{L1,L2} has proved essentially best possible bounds for convex hypersurfaces (codimension 1). Later, Beresnevich, Vaughan, Velani, Zorin \cite{BVVZ} and Simmons \cite{Si}  have made further progress for submanifolds satisfying various rank or curvature conditions, which only hold generically for submanifolds with  sufficiently large dimensions in terms of the ambient dimension. Recently, J. Liu and the author \cite{HL1} have obtained optimal  results for affine subspaces of $\R^n$  satisfying certain diophantine type conditions. 

In this paper, we are primarily concerned with the above counting problem for space curves in $\R^3$, which has been under people's radar for quite a while, but seems to have resisted all the attacks so far. Here we make some progress on this problem.

By the inverse function theorem, any $C^3$ space curve may be parametrized in the Monge form locally. If the curve is also compact, it can be covered by finitely many subcurves, each presented in the Monge form. Hence bounds for the global counting problem can be obtained from local ones, at the expense of losing a constant factor. Therefore, after proper dilation and translation, we may assume that the curve $\mathcal{C}$ is parametrized by
\begin{equation*}
\{(x, f_1(x), f_2(x)):x\in [0,1]\},\quad  \textrm{where } f_1, f_2\in C^3([0,1]),
\end{equation*}
and moreover the torsion of $\mathcal{C}$ is nowhere vanishing if and only if
\begin{equation}\label{e1}
\begin{vmatrix}
1 & f_1'(x) & f_2'(x) \\ 
0 & f_1''(x) & f_2''(x) \\ 
0 & f_1'''(x) & f_2'''(x) 
\end{vmatrix}
=
\begin{vmatrix}
f_1''(x) & f_2''(x) \\ 
f_1'''(x) & f_2'''(x)
\end{vmatrix}\not=0, \quad\textrm{for all }x\in [0,1].
\end{equation}

Let $$A(q,\delta):=\#\{a\in [0,q]\cap\mathbb{Z}:\|qf_1(a/q)\|<\delta, \|qf_2(a/q)\|<\delta\}.$$ Roughly speaking, $A(q,\delta)$ counts the number of rational points with denominator $q$ lying in the $O(\delta/q)$ neighborhood of $\mathcal{C}$, or equivalently this is the same as counting lattice points within distance $O(\delta)$ to the dilation $q\mathcal{C}$ of the curve $\mathcal{C}$. A simple probabilistic heuristic shows that one expects to have $A(q,\delta)\asymp \delta^2q$, which of course breaks down when $\delta\to0$ and $q$ is fixed. The following theorem confirms this heuristic when $\delta q^{\frac15}(\log q)^\frac25\to\infty$.

\begin{thm}\label{t1}
Let $\mathcal{C}$ be a compact $C^3$ curve in $\mathbb{R}^3$ with nonvanishing torsion. Then
for any $\delta\in(0,1/2)$ and $q\ge1$, we have
$$
A(q,\delta)\ll_{\mathcal{C}}\delta^2q+q^{\frac35}(\log q)^{\frac45}.
$$
Furthermore, there exist positive constants $C$ and $Q_0$ such that when $\delta\ge Cq^{-\frac15}(\log q)^{\frac25}$ and $q\ge Q_0$ we have
$$
A(q,\delta)\gg_{\mathcal{C}} \delta^2q.
$$
\end{thm}

Our Theorem \ref{t1} represents the first nontrivial result of its kind about curves in $\R^3$. The term $\delta^2 q$ is the heuristic main term hence cannot be dispensed with. It remains to be seen whether the other term $q^{\frac35}(\log q)^{\frac45}$  is subject to improvement. It is not unlikely that 
$$
A(q,\delta)\ll\delta^2q+q^{\frac12+\varepsilon}
$$
or even 
\begin{equation}\label{e3}
A(q,\delta)\ll\delta^2q+q^{\frac13+\varepsilon}.
\end{equation}
We may as well expect that the lower bound $A(q,\delta)\gg\delta^2q$ holds provided that $\delta\gg q^{-\frac14+\varepsilon}$ (or $\delta\gg q^{-\frac13+\varepsilon}$).

In view of the cubic Veronese curve $\mathcal{V}_3=\{(x,x^2,x^3):x\in[0,1]\}$, we must have the lower bound 
$$
\#(q\mathcal{V}_3\cap\mathbb{Z}^3)\gg q^{\frac13}.
$$
Therefore the conjectural bound \eqref{e3} is the best that one can hope for. Nevertheless, further improvement should be possible if one averages over $q$. Indeed, it is reasonable to conjecture that
$$
\sum_{q\le Q} A(q,\delta)\ll\delta^2Q^2+Q^{1+\varepsilon}
$$
or even 
$$
\sum_{q\le Q} A(q,\delta)\ll\delta^2Q^2+Q^{\frac23+\varepsilon}.
$$
The above conjectures, if true, would have significant consequences in metric diophantine approximations on nondegenerate curves in $\R^3$. We only remark in passing that estimating the number of lattice/rational points close to a manifold is a very active area of research and refer the interested readers to the papers \cite{bere, BDV, BVVZ, BVVZ1, BZ, hua1, hua3, hux, Le, Si, VV} for more background and recent developments. 

An alert reader may wonder whether the torsion condition \eqref{e1} is necessary for Theorem \ref{t1} to hold true. Here we provide an example to answer the above question affirmatively. Consider the embedded parabola $\{(x,x^2,0)|x\in[0,1]\}$ in $\R^3$, which clearly fails the torsion condition \eqref{e1}. In this case, the second inequality in the definition of $A(q,\delta)$ is always true, hence the counting problem is equivalent to its analogue in $\R^2$ for the standard parabola. The latter problem has been studied recently by H. Li and the author \cite[Theorem 2]{HL2}\footnote{The result there is more precise than the one quoted here.}, and we have
$$
A(q,\delta)=2\delta q+O(q^{\frac12+\varepsilon}),\quad \text{for any }\varepsilon>0,
$$
which is clearly incompatible with the upper bound in Theorem \ref{t1}. This shows that our Theorem \ref{t1} cannot hold for this embedded parabola, and therefore the nonvanishing torsion condition assumed in Theorem \ref{t1} cannot be completely dispensed with.

The novel feature of our approach is an induction scheme which enables us to reduce the lattice points counting problem for space curves to one that is relevant to counting rational points near planar curves. More precisely, the starting point of our proof of Theorem \ref{t1} is based on an argument of Sprind\v{z}uk \cite[\S 2.9]{Sp}, which is also revisited in \cite{BVVZ, Si}. In a nutshell, the idea is that we approximate the curve by some short line segments in such a way that a thin neighborhood of the curve more or less coincides with a thin neighborhood of the broken line segments. This way, the counting problem is locally linearized, therefore can be better handled by analytic techniques. However, in order to add up the contributions from those linear patches and obtain enough savings to claim victory, one has to assume some rank or curvature conditions in Sprind\v{z}uk's original argument and its later variations, which unfortunately completely exclude curves in $\mathbb{R}^n$ with $n\ge3$. It is at this stage that we deviate from all previous approaches and have to reply on a very delicate analysis which utilizes information from one lower dimension. Finally, we use essentially the best possible result in $\mathbb{R}^2$ c.f. Proposition \ref{l5} as the initial input. In principle, one may use this induction scheme to obtain bounds for lattice points close to  a curve in $\mathbb{R}^n$ with $n\ge4$. Unfortunately, the quality of the estimates rendered this way deteriorates fairly fast as $n$ increases. Therefore, we decide not to pursue the full potential of our method herein, and only focus on the pivotal case $n=3$.

\section{Preliminary Lemmata}

The main purpose of this section is to prove the Proposition \ref{l5} below, whose argument draws upon our earlier work \cite{hua2}.

Let $I=[\xi,\eta]$. Suppose that $f\in C^2(I)$ satisfies
\begin{equation}\label{e4}
0<c_1\le|f''(x)|\le c_2.
\end{equation}
Let
$$
\mu(j_1, j_2, \lambda)=\{x\in I: \|j_1x+j_2f(x)\|<\lambda\}
$$
and 
\begin{equation*}
\mu(j_1, j_2, p, \lambda)=\{x\in I: |j_1x+j_2f(x)-p|<\lambda\}.
\end{equation*}

\begin{prop}\label{l5} For positive integer $J$ and $\lambda\in(0,\frac12)$, 
we have
$$
\sum_{\substack{|j_1|,|j_2|\le J\\(j_1,j_2)\not=(0,0)}}|\mu(j_1, j_2, \lambda)|\ll \lambda J^2+\lambda^{\frac12}J^{\frac12}\log J.
$$
\end{prop}

\subsection{Lemmata} We state some lemmata first, which will be used in the proof of Proposition \ref{l5}.

\begin{lem}[{\cite[Lemma 9.7]{Ha}}]\label{l2}
Let $h(x)\in C^2(I)$ be such that $\min_{x\in I}|h'(x)|=\delta_1$ and $\min_{x\in I}|h''(x)|=\delta_2$. For $\tau>0$, define
$$E(\tau):=\{x\in I:|h(x)|<\tau\}.$$ Then we have
$$|E(\tau)|\ll\min\left(\frac\tau{\delta_1},\sqrt{\frac{\tau}{\delta_2}}\,\right).$$
\end{lem}


\begin{lem}[{\cite[Lemma 4]{hua2}}]\label{l3}
Suppose that $\phi$ has a continuous second derivative on a bounded interval $K$ which is bounded away from 0, and let $\delta\in(0,\frac14)$. Then for any $\varepsilon>0$ and $U\ge1$,
$$\sum_{U\le u<2U}\sum_{\substack{t/u\in K\\\|u\phi(t/u)\|<\delta}}1\ll_{\varepsilon, K} \delta^{1-\varepsilon}U^2+U\log(2U).$$ 
\end{lem}

\begin{lem}[{\cite[Lemma 5]{hua2}}]\label{l4}
Under the same conditions with Lemma \ref{l3}, for $\Lambda\in(0,1)$ and $U\ge1$,
$$\sum_{U\le u<2U}\sum_{\substack{t/u\in K\\\|u\phi(t/u)\|\ge\delta}}\left\|u\phi\left(\frac{t}{u}\right)\right\|^{-\Lambda}\ll_K U^2+\delta^{-\Lambda}U\log(2U).$$ 
\end{lem}

\subsection{Proof of Poposition \ref{l5}}

We will divide the set $[-J,J]\cap\mathbb{Z}^2\backslash(0,0)$ of all possible choices for $(j_1,j_2)$ into a couple of subsets. We also notice that for given $(j_1,j_2)$, there are only finitely many $p\in\mathbb{Z}$ such that $\mu(j_1,j_2,p,\lambda)\not=\emptyset$. Indeed, such $p$ must satisfy
\begin{equation}\label{e3.3}
|p|\le C \max\{|j_1|,|j_2|\}
\end{equation}
where $$C=\max_{x\in I}\{|x|+|f(x)|+1\}.$$

Let
$$M=1+\max_{x\in I}|f'(x)|$$
and
$$
\Theta=[-J,J]^2\cap\mathbb{Z}\times\mathbb{Z}\setminus(0,0),
$$
$$\Theta_1=\{(j_1,j_2)\in \Theta:|j_1|>2M|j_2|\},$$
$$\Theta_2=\Theta\setminus\Theta_1.$$

We consider the case $(j_1,j_2)\in\Theta_1$ first. In this case, we have
$$|j_1+j_2f'(x)|\ge|j_1|-M|j_2|\ge\frac{|j_1|}2.$$
Now by Lemma \ref{l2}, we know
$$
|\mu(j_1,j_2,p,\lambda)|\ll \frac{\lambda}{|j_1|}.
$$
Moreover, for a given $j_1$, there are at most $\ll j_1^2$ possible choices for $j_2$ and $p$. 
Therefore 
$$
\sum_{\substack{(j_1,j_2)\in\Theta_1\\p\in\mathbb{Z}}}|\mu(j_1,j_2,p,\lambda)|
{\ll}\lambda J^{2}.
$$

Now for $(j_1,j_2)\in\Theta_2$, clearly  we have $j_2\not=0$ and $|j_1|\le 2M |j_2|$. For convenience, we may extend the definition of $f(x)$ to $\mathbb{R}$  by taking the second order Taylor expansions at the end points of $I$. Namely
let 
$$f(x)=f(\eta)+f'(\eta)(x-\eta)+\frac{f''(\eta)}2(x-\eta)^2$$
when $x>\eta$
and
$$f(x)=f(\xi)+f'(\xi)(x-\xi)+\frac{f''(\xi)}2(x-\xi)^2$$
when $x<\xi$.
Clearly the extended $f$ satisfies $f\in C^2(\mathbb{R})$ and $c_1\le |f''|\le c_2$. Since $f''$ does not change sign throughout $\mathbb{R}$, $f'$ is strictly monotonic on $\mathbb{R}$ and has range $(-\infty,\infty)$. Let $g(y):\mathbb{R}\rightarrow\mathbb{R}$ be the inverse function of $-f'(x)$. To this end,
let $$x_0:=g(j_1/j_2)$$ 
which is the unique point $x_0\in\mathbb{R}$ such that
\begin{equation}\label{e3.4}
j_1+j_2f'(x_0)=0.
\end{equation}
Let $K=[-2M,2M]$ and $I'=g(K)\supset I$. So $x_0\in I'$. 
Note that
$$g'(y)=\frac{-1}{f''(g(y))}$$
and hence that
\begin{equation}\label{e3.12}
c_2^{-1}\le|g'(y)|\le c_1^{-1}
\end{equation}
for all $y\in\mathbb{R}$. Thus by the mean value theorem
\begin{equation}\label{e3.5}
|I'|\le c_1^{-1}|K|= 4c_1^{-1}M.
\end{equation}

Now, we let
$$F(x)=j_1x+j_2f(x)$$ with $j_1/j_2\in K$. Then by \cite[Lemma 3]{hua2} we have the following lemma.
\begin{lem}\label{l1}
$$|F'(x)|\asymp |j_2(F(x)-F(x_0))|^{1/2}$$
where the $\asymp$ constants depend only on $c_1$, $c_2$.
\end{lem}

Now let $p_0$ be the unique integer such that
\begin{equation}\label{e3.8}
-\frac12<F(x_0)-p_0\le\frac12.
\end{equation}
If $p\not=p_0$, then for $x\in\mu(j_1,j_2,p,\lambda)$
\begin{align*}
|F(x)-F(x_0)|&=|p-p_0+F(x)-p-F(x_0)+p_0|\\
&\ge|p-p_0|-|F(x)-p|-|F(x_0)-p_0|\\
&\overset{\eqref{e3.8}}{\ge}|p-p_0|-\lambda-1/2\\
&\ge\frac13|p-p_0|
\end{align*}
provided that 
\begin{equation}\label{e3.9}
\lambda\le \frac1{8}.
\end{equation} 
Proposition \ref{l5} is trivial when $\lambda>\frac18$, so without loss of generality we may assume \eqref{e3.9} holds.

By Lemma \ref{l1} and then Lemma \ref{l2}, we get, when $p\not=p_0$ 
\begin{equation}\label{e10}
|\mu(j_1,j_2,p,\lambda)|\ll \lambda(|j_2||p-p_0|)^{-1/2}.
\end{equation}
Therefore for fixed $(j_1,j_2)\in\Theta_2$ 
\begin{align*}
\sum_{p\not=p_0}|\mu(j_1,j_2,p,\lambda)|&\ll \sum_{p\not=p_0}\lambda(j_2|p-p_0|)^{-1/2}\\
&\overset{\eqref{e3.3}}{\ll}\lambda j_2^{-1/2}j_2^{1-1/2}\\
&\ll \lambda.
\end{align*}
 Hence
$$\sum_{\substack{(j_1,j_2)\in\Theta_2\\p\not=p_0}}|\mu(j_1,j_2,p,\lambda)|\ll \lambda J^{2}.$$

Hitherto, we are left with the most difficult case $p=p_0$. For $x\in\mu(j_1,j_2,p_0,\lambda)$, we have
\begin{align*}
|F(x)-F(x_0)|&=|F(x)-p_0+p_0-F(x_0)| \\
&\ge\|F(x_0)\|-\lambda\\
&\ge\frac12\|F(x_0)\|
\end{align*}
provided that
$$\|F(x_0)\|\ge2\lambda.$$
By Lemma \ref{l1} and the inequality above
$$|F'(x)|\gg (j_2\|F(x_0)\|)^{1/2}.$$
Then by Lemma \ref{l2}
\begin{equation}\label{e3.10}
|\mu(j_1,j_2,p_0,\lambda)|\ll\frac{\lambda}{j_2^{1/2}}\|F(x_0)\|^{-1/2},
\end{equation}
when $\|F(x_0)\|\ge2\lambda$. 

In the case if $\|F(x_0)\|<2\lambda$, we simply use the upper bound
\begin{equation}\label{e3.11}
|\mu(j_1,j_2,p_0,\lambda)|\overset{\text{Lem.}\ref{l2}}{\ll} \sqrt{\frac{\lambda}{j_2}}.
\end{equation}

We now define the dual curve $f^*(y)$ of $f(x)$, whose derivative is the inverse function of $-f'(x)$. Namely
$$f^*(y):=yg(y)+f(g(y)).$$
It is readily seen that
\begin{equation}\label{e3.13}
(f^*)'(y)=g(y)+yg'(y)+f'(g(y))g'(y)=g(y)
\end{equation}
and that
\begin{equation}\label{e3.14}
j_2f^*(j_1/j_2)=j_1g(j_1/j_2)+j_2f(g(j_1/j_2))=j_1x_0+j_2f(x_0)=F(x_0).
\end{equation}

Finally, we are poised to treat the sum
\begin{equation}\label{e3.15}
\sum_{(j_1,j_2)\in\Theta_2}|\mu(j_1,j_2,p_0,\lambda)|.
\end{equation}
 We further divide this into two cases, namely $\|F(x_0)\|<2\lambda$ and $\|F(x_0)\|\ge2\lambda$.

For any nonnegative integer $k$,
\begin{align*}
&\sum_{\substack{(j_1,j_2)\in\Theta_2\\2^{k}\le |j_2|<2^{k+1}\\\|F(x_0)\|<2\lambda}}|\mu(j_1,j_2,p_0,\lambda)|\\
\overset{\eqref{e3.11}\&\eqref{e3.14}}{\ll}& \sum_{2^{k}\le j_2<2^{k+1}}\sum_{\substack{j_1/j_2\in K\\\|j_2f^*(j_1/j_2)\|<2\lambda}}\left(\frac{\lambda}{2^k}\right)^{\frac{1}2}\\
\overset{\text{Lem.}\ref{l3}}{\ll}&\Big(\lambda^{1-\varepsilon}2^{2k}+(k+1)2^k\Big)\left(\frac{\lambda}{2^k}\right)^{\frac{1}2}\\
\ll&\lambda^{\frac{3}2-\varepsilon}2^{3k/2}+(k+1)\lambda^\frac{1}2 2^{k/2}.
\end{align*}
By summing over $0\le k\le \log_2 J$, we have
\begin{equation}\label{e3.17}
\sum_{\substack{(j_1,j_2)\in\Theta_2\\\|F(x_0)\|<2\lambda}}|\mu(j_1,j_2,p_0,\lambda)|\ll \lambda^{\frac32-\varepsilon}J^{\frac32}+\lambda^{\frac12}J^{\frac12}\log J.
\end{equation}

The other case can be treated in a similar fashion. 
\begin{align*}
&\sum_{\substack{(j_1,j_2)\in\Theta_2\\2^{k}\le |j_2|<2^{k+1}\\\|F(x_0)\|\ge2\lambda}}|\mu(j_1,j_2,p_0,\lambda)|\\
\overset{\eqref{e3.10}\&\eqref{e3.14}}{\ll}& \sum_{2^{k}\le j_2<2^{k+1}}\sum_{\substack{j_1/j_2\in K\\\|j_2f^*(j_1/j_2)\|\ge2\lambda}}\left(\frac{\lambda}{2^{k/2}}\right)\|j_2f^*(j_1/j_2)\|^{-1/2}\\
\overset{\text{Lem.}\ref{l4}}{\ll}&\left(\lambda^{-1/2}2^{k}(k+1)+2^{2k}\right)\frac{\lambda}{2^{k/2}}\\
\ll&\lambda^{1/2}2^{k/2}(k+1)+\lambda 2^{3k/2}.
\end{align*}
Again by summing over $0\le k\le \log_2 J$, we obtain
\begin{equation}\label{e3.18}
\sum_{\substack{(j_1,j_2)\in\Theta_2\\\|F(x_0)\|\ge2\lambda}}|\mu(j_1,j_2,p_0,\lambda)|\ll \lambda^{\frac12}J^{\frac12}\log J+\lambda J^{\frac32}.
\end{equation}

\section{Proof of Theorem \ref{t1}}

Let
$$
\mathcal{A}(q,\delta)=\{a\in [0,q]\cap\mathbb{Z}:\|qf_1(a/q)\|<\delta, \|qf_2(a/q)\|<\delta\}.
$$

In view of the torsion condition \eqref{e1}, either $f_1''$ or $f_2''$ must be bounded away from 0 in a sufficiently small neighborhood of any $x\in [0,1]$. Therefore, without loss of generality, we may prove Theorem \ref{t1} under the additional assumption that $|f_1''|\ge c_3>0$. The general case would follow immediately by compactness.

Let $c_4=\max(\|f_1\|_{C^3([0,1])},\|f_2\|_{C^3([0,1])},1)$, $q_0=\lfloor(2^{-1}c_4^{-1}\delta q)^{\frac12}\rfloor$ and $r=\lfloor q/q_0\rfloor$. For the time being, we suppose $\delta\ge 2c_4q^{-1}$ so that $q_0\ge1$. For each $a\in [0,q]\cap\mathbb{Z}$ write $a=q_0s+a_0$ with $a_0=a_0(a)\in[0, q_0)\cap\mathbb{Z}$ and $s=s(a)\in[0, r]\cap\mathbb{Z}$.
Let
$$
\mathcal{A}(q,\delta,s)=\{a\in\mathcal{A}(q,\delta): s(a)=s\}
$$
and
$$
{A}(q,\delta,s)=\#\mathcal{A}(q,\delta,s).
$$
By the Taylor theorem, when $a\in\mathcal{A}(q,\delta,s)$, we have for $i\in\{1,2\}$
\begin{equation*}
\left|f_i\left(\frac{a}q\right)-f_i\left(\frac{q_0s}q\right)-\frac{a_0}qf_i'\left(\frac{q_0s}q\right)\right|\le \frac{c_4}2\cdot\frac{a_0^2}{q^2}\le \frac{c_4q_0^2}{2q^2}< \frac{\delta}{2q}.
\end{equation*}
Then it follows by the triangle inequality that for $i\in\{1,2\}$
\begin{equation}\label{e4.1}
\Bigg|\left\|qf_i\left(\frac{q_0s}q\right)+{a_0}f_i'\left(\frac{q_0s}q\right)\right\|-\left\|qf_i\left(\frac{a}q\right)\right\|\Bigg|\le\frac{\delta}2.
\end{equation}
Let $$\mathcal{B}(q,\delta,s):=\left\{a_0\in[0,q_0)\cap\mathbb{Z}:\left\|qf_i\left(\frac{q_0s}q\right)+{a_0}f_i'\left(\frac{q_0s}q\right)\right\|<\delta,\quad 1\le i\le 2\right\},$$
$$
B(q,\delta,s):=\#\mathcal{B}(q,\delta,s),
$$
$$
B_1(q,\delta):=\sum_{0\le s\le r}B(q,\delta,s),
$$
and
$$
B_2(q,\delta):=\sum_{0\le s< r}B(q,\delta,s).
$$

We then observe from \eqref{e4.1} that 
$$
a=q_0s+a_0\in\mathcal{A}(q,\delta,s)\Rightarrow a_0\in\mathcal{B}(q,3\delta/2,s),\quad\textrm{when } s\le r
$$
and
$$
 a_0\in\mathcal{B}(q,\delta/2,s)\Rightarrow a=q_0s+a_0\in\mathcal{A}(q,\delta,s)\quad\textrm{when } s<r.
$$
Therefore we obtain
$$
{A}(q,\delta,s)\le{B}(q,3\delta/2,s),\quad\textrm{when } s\le r
$$
and
$$
 {B}(q,\delta/2,s)\le {A}(q,\delta,s)\quad\textrm{when } s<r,
$$
and hence
\begin{equation}\label{e4.8}
B_2(q,\delta/2)\le A(q,\delta)\le B_1(q,3\delta/2).
\end{equation}

Now to estimate ${B}(q,\delta,s)$, we recall some basic properties of Selberg's magic functions. See \cite[Chapter 1]{mo} for details about the construction of these functions.

Let $\Delta=(\alpha, \beta)$ be an arc of $\R/\Z$ with $\alpha<\beta<\alpha+1$,  and $\chi_{_\Delta}(x)$ be its  characteristic function. Then there exist finite trigonometric polynomials of degree at most $J$
$$S^{\pm}_J(x)=\sum_{|j|\le J}b_j^{\pm}e(jx)$$ such that 
$$
S^{-}_J(x)\le \chi_{_\Delta}(x)\le S^{+}_J(x)
$$
and
$$b_0^\pm=\beta-\alpha\pm\frac1{J+1}
$$
and
$$|b_j^\pm|\le\frac1{J+1}+\min\left(\beta-\alpha,\frac1{\pi|j|}\right)
$$
for $0<|j|\le J$. Here $e(x)=e^{2\pi ix}$.

Here we choose $\alpha=-\delta$, $\beta=\delta$ and $J=\lfloor\frac1{\delta}\rfloor$. With such choices, observe that $|b_j^\pm|\le3\delta$ for all $j\in\mathbb{Z}$ and $b_0^-\ge\delta$.  

Let $$
F_i(s,a_0):=qf_i\left(\frac{q_0s}q\right)+{a_0}f_i'\left(\frac{q_0s}q\right),\quad 1\le i\le 2.$$
Then
\begin{align*}
B(q,\delta,s)=& \sum_{0\le a_0<q_0}\chi_{_\Delta}(F_1(s,a_0))\chi_{_\Delta}(F_2(s,a_0))\\
\le&\sum_{0\le a_0<q_0}S^+_J(F_1(s,a_0))S^+_J(F_2(s,a_0))\\
=&\sum_{|j_1|,|j_2|\le J}b_{j_1}^+b_{j_2}^+\sum_{0\le a_0<q_0}e\left(\sum_{i=1}^2j_iF_i(s,a_0)\right)\\
\le&9\delta^2 q_0+\sum_{\substack{|j_1|,|j_2|\le J\\(j_1,j_2)\not=(0,0)}}b_{j_1}^+b_{j_2}^+\sum_{0\le a_0<q_0}e\left(\sum_{i=1}^2j_iF_i(s,a_0)\right),
\end{align*}
where in the last inequality we single out the zero frequency term $j_1=j_2=0$ which gives the heuristic main term for $B(q,\delta,s)$.
Similarly
\begin{align*}
B(q,\delta,s)=& \sum_{0\le a_0<q_0}\chi_{_\Delta}(F_1(s,a_0))\chi_{_\Delta}(F_2(s,a_0))\\
\ge&\sum_{0\le a_0<q_0}S^-_J(F_1(s,a_0))S^-_J(F_2(s,a_0))\\
=&\sum_{|j_1|,|j_2|\le J}b_{j_1}^-b_{j_2}^-\sum_{0\le a_0<q_0}e\left(\sum_{i=1}^2j_iF_i(s,a_0)\right)\\
\ge&\delta^2 q_0+\sum_{\substack{|j_1|,|j_2|\le J\\(j_1,j_2)\not=(0,0)}}b_{j_1}^-b_{j_2}^-\sum_{0\le a_0<q_0}e\left(\sum_{i=1}^2j_iF_i(s,a_0)\right).
\end{align*}


Note that
\begin{align*}
\left|\sum_{0\le a_0<q_0}e\left(\sum_{i=1}^2j_iF_i(s,a_0)\right)\right|=&\left|\sum_{0\le a_0<q_0}e\left(\sum_{i=1}^2j_i{a_0}f_i'\left(\frac{q_0s}q\right)\right)\right|\\
\le&\min\left(q_0,\left\|\sum_{i=1}^2j_if_i'\left(\frac{q_0s}q\right)\right\|^{-1}\right),
\end{align*}
where in the last inequality we use the well known linear exponential sum estimate
$$
\left|\sum_{n\le N}e(n\gamma)\right|\le\min(N, \|\gamma\|^{-1}).
$$

Thus
\begin{equation}\label{e4.2}
B_1(q,\delta)\le9\delta^2(q+q_0)+E(q,\delta)
\end{equation}
and
\begin{equation}\label{e4.3}
B_2(q,\delta)\ge\delta^2 (q-q_0)-E(q,\delta),
\end{equation}
where
$$
E(q,\delta)=9\delta^2\sum_{0\le s\le r}\sum_{\substack{|j_1|,|j_2|\le J\\(j_1,j_2)\not=(0,0)}}\min\left(q_0,\left\|\sum_{i=1}^2j_if_i'\left(\frac{q_0s}q\right)\right\|^{-1}\right).$$

\par

For a given integer $s$, consider the intervals $I_s=[s-1/2,s+1/2]$, unless $s=0$ or $s=r$ in which case we consider $[s, s+1/2]$ or $[s-1/2,s]$ respectively. For $\alpha\in I_s$ we have
$$
\left|f_i'\left(\frac{q_0s}q\right)-f_i'\left(\frac{q_0\alpha}q\right)\right|\le c_4\frac{q_0}{2q},\quad 1\le i\le 2.
$$ 
Hence
$$
\left|\sum_{i=1}^2j_i\left(f_i'\left(\frac{q_0s}q\right)-f_i'\left(\frac{q_0\alpha}q\right)\right)\right|\le c_4J\frac{q_0^2}{qq_0}\le c_4\frac1{\delta}\frac{\delta q}{2c_4qq_0}\le \frac1{2q_0}.
$$
Thus
\begin{equation}\label{e4.4}
\min\left(q_0,\left\|\sum_{i=1}^2j_if_i'\left(\frac{q_0s}q\right)\right\|^{-1}\right)\le2 \min\left(q_0,\left\|\sum_{i=1}^2j_if_i'\left(\frac{q_0\alpha}q\right)\right\|^{-1}\right).
\end{equation}
Now we integrate over $\alpha\in I_s$ and then sum over $s\in[0, r]$. This way we obtain
\begin{align}
&\nonumber\sum_{0\le s\le r}\min\left(q_0,\left\|\sum_{i=1}^2j_if_i'\left(\frac{q_0s}q\right)\right\|^{-1}\right)\\\nonumber\le&\left(2\int_0^{\frac12}+\int_{\frac12}^{r-\frac12}+2\int_{r-\frac12}^r \right)2\min\left(q_0,\left\|\sum_{i=1}^2j_if_i'\left(\frac{q_0\alpha}q\right)\right\|^{-1}\right)d\alpha\\\nonumber\le&4\int_0^r \min\left(q_0,\left\|\sum_{i=1}^2j_if_i'\left(\frac{q_0\alpha}q\right)\right\|^{-1}\right)d\alpha\\\nonumber
\le&4\int_0^{\frac{q}{q_0}} \min\left(q_0,\left\|\sum_{i=1}^2j_if_i'\left(\frac{q_0\alpha}q\right)\right\|^{-1}\right)d\alpha\\\label{e4.5}
\ll&\frac{q}{q_0}\int_{I} \min\left(q_0,\left\|j_1\beta+j_2f(\beta)\right\|^{-1}\right)d\beta
\end{align}
where
$$
I=[\inf f_1',\sup f_1']
$$
and
$$
f=f_2'\circ(f_1')^{-1}.
$$
In the last inequality, we make the change of variable $\beta=f'_1(q_0\alpha/q)$ and use the assumption that $f''_1$ is bounded away from 0. 
\par
Now we verify that the torsion condition \eqref{e1} implies that $f''$ is bounded away from zero. Since $f\circ f_1'=f_2'$, we have
$$
(f'\circ f_1')\cdot f_1''=f_2'',
$$ 
$$
(f''\circ f_1')\cdot (f_1'')^2+(f'\circ f_1')\cdot f_1'''=f_2'''
$$
and hence
$$
f''\circ f_1'=\frac{f_1''f_2'''-f_2''f_1'''}{(f_1'')^3}
$$
which is bounded away from 0 in view of \eqref{e1}.

We are ready to estimate the resulting integral from \eqref{e4.5}. Our strategy is to decompose $I$ into subsets defined by the inequalities $\left\|j_1\beta+j_2f(\beta)\right\|<\frac{1}{q_0}$ and $\frac{2^{k-1}}{q_0}\le\left\|j_1\beta+j_2f(\beta)\right\|<\frac{2^k}{q_0}$ for $1\le k\le \log_2q_0$. This naturally leads to the set $\mu(j_1,j_2,\lambda)$ defined in the previous section, for which we may invoke Proposition \ref{l5}. 
Therefore, we have

\begin{align}
&\nonumber\sum_{\substack{|j_1|,|j_2|\le J\\(j_1,j_2)\not=(0,0)}}\int_{I} \min\left(q_0,\left\|j_1\beta+j_2f(\beta)\right\|^{-1}\right)d\beta\\
&\nonumber\ll\sum_{\substack{|j_1|,|j_2|\le J\\(j_1,j_2)\not=(0,0)}}\left(q_0|\mu(j_1,j_2,1/q_0)|+\sum_{1\le k\le\log_2q_0}\int_{\frac{2^{k-1}}{q_0}\le\left\|j_1\beta+j_2f(\beta)\right\|<\frac{2^k}{q_0}} \left\|j_1\beta+j_2f(\beta)\right\|^{-1}d\beta\right)\\
&\nonumber\ll \sum_{\substack{|j_1|,|j_2|\le J\\(j_1,j_2)\not=(0,0)}}\sum_{0\le k\le \log_2 q_0}\frac{q_0}{2^k}|\mu(j_1, j_2, 2^k/q_0)|\\
&\nonumber\overset{\textrm{Prop. }\ref{l5}}{\ll}\sum_{0\le k\le \log_2 q_0}q_02^{-k}(2^kq_0^{-1}J^2+2^{k/2}q_0^{-1/2}J^{1/2}\log J)\\
&\ll J^2\log q_0+q_0^{1/2}J^{1/2}\log J.\label{e4.6}
\end{align}
Hence we conclude from \eqref{e4.4}, \eqref{e4.5} and \eqref{e4.6},  that
\begin{align*}
E(q,\delta)
\le& C_1\delta^2\frac{q}{q_0}(J^2\log q_0+q_0^{1/2}J^{1/2}\log J)\\
\le &C_2(\delta^{-1/2}q^{1/2}\log q+\delta^{5/4}q^{3/4}\log \frac1\delta),
\end{align*}
where the second inequality follows on noting that $J=\lfloor 1/\delta\rfloor$ and $q_0=\lfloor(2^{-1}c_4^{-1}\delta q)^{\frac12}\rfloor$.

Hence we obtain from \eqref{e4.2} that 
$$
B_1(q,\delta)\ll \delta^2 q+\delta^{-1/2}q^{1/2}\log q+\delta^{5/4}q^{3/4}\log \frac1\delta,\quad \text{when } \delta\ge 2c_4q^{-1}.
$$
Now observe that for fixed $q$, $B_1(q,\delta)$ is increasing in $\delta$. Let $\delta_0=q^{-\frac15}(\log q)^{\frac25}$. Therefore when $\delta<2c_4\delta_0$,
$$
B_1(q,\delta)\le B_1(q,2c_4\delta_0)\ll q^{\frac35}(\log q)^{\frac45},
$$
and when $\delta\ge 2c_4\delta_0$,
$$
B_1(q,\delta)\ll \delta^2q.
$$
In any case, we have for all $\delta\in(0,1/2)$ that
\begin{equation}\label{e4.9}
B_1(q,\delta)\ll \delta^2q+q^{\frac35}(\log q)^{\frac45}.
\end{equation}

\par 
On the other hand, there exist positive constants $C, Q_0$ such that when $\delta\ge C\delta_0$ and $q\ge Q_0$ we have
$$
\delta^2q_0+C_2(\delta^{-1/2}q^{1/2}\log q+\delta^{5/4}q^{3/4}\log \frac1\delta)
\le  \frac1{100}\delta^2q.
$$
Therefore, we obtain the lower bound from \eqref{e4.3},
\begin{equation}\label{e4.10}
B_2(q,\delta)\ge 0.99 \delta^2q, \quad\text{when }\delta\ge C\delta_0, q\ge Q_0.
\end{equation}

Now the proof follows by combining  \eqref{e4.8}, \eqref{e4.9} and \eqref{e4.10}.

\proof[Acknowledgments]
The author is grateful to the anonymous referee for carefully reading the manuscript and providing helpful suggestions to improve the presentation.

\end{document}